\documentclass{amsart}
\usepackage{hyperref}
\usepackage{graphicx,amsmath,amssymb,amsthm}
\usepackage{color}

\def\<#1,#2>{\langle#1,#2\rangle}
\newcommand{\appr}{^{^{_\sim}}}

\newcommand{\mrm}[1]{\text{\rm #1}}
\newcommand{\dd}{\,d}
\newcommand{\col}{\operatorname{col}}

\newcommand{\new}[1]{{\em #1}}

\newcommand{\R}{\mathbb{R}}
\newcommand{\rbar}{\overline{\R}}

\newcommand{\sK}{\mathcal{K}}

\newcommand{\im}{\mathrm{im}\,}

\newcommand{\projimker}[2]{\Pi_{#1}^{#2}}
\newcommand{\rmax}{\mathbb{R}_{\max}}
\newcommand{\rmaxb}{\overline{\R}_{\max}}
\newcommand{\rminb}{\overline{\R}_{\min}}

\newcommand{\supp}{\mathop{\text{\Large$\vee$}}}
\newcommand{\inff}{\mathop{\text{\Large$\wedge$}}}
\newcommand{\maxx}{\max}

\newcommand{\comp}{\circ}

\newcommand{\bydef}{\stackrel{\mathrm{def}}{=}}
\newcommand{\set}[2]{\{#1\mid\,#2\}}
\newcommand{\lres}{\backslash}

\newcommand{\sh}{^{\sharp}}

\newcommand{\calu}{\mathcal{U}}
\newcommand{\calv}{\mathcal{V}}
\newcommand{\calw}{\mathcal{W}}
\newcommand{\calx}{\mathcal{X}}
\newcommand{\caly}{\mathcal{Y}}
\newcommand{\calz}{\mathcal{Z}}
\def\<#1,#2>{\langle#1\mid #2\rangle}
\newtheorem{thm}{Theorem}
\newtheorem{prop}[thm]{Proposition}
\newtheorem{lem}[thm]{Lemma}

\newtheorem{exmp}[thm]{Example}
\newtheorem{rem}[thm]{Remark}
\newtheorem{defin}[thm]{Definition}
\title[The max-plus finite element method]{The max-plus finite element method for optimal control problems: further approximation results}
\author{Marianne Akian}
\author{St\'ephane Gaubert}
\author{Asma Lakhoua}
\address{INRIA, Domaine de Voluceau, 78153 Le Chesnay C\'edex, France}
\email{\{Marianne.Akian,Stephane.Gaubert,Asma.Lakhoua\}@inria.fr}
\date{September 12, 2005. Prepared for CDC-ECC'05.}
\keywords{Max-plus algebra, tropical semiring, Hamilton-Jacobi equation,
weak formulation, residuation, projection, idempotent semimodules, finite element method.}

\subjclass[2000]{Primary 49L20; Secondary 65M60, 06A15, 12K10}
\begin{document}
\begin{abstract}
We develop the max-plus finite element method to solve finite horizon deterministic optimal control problems. This method, that we introduced in a previous work, relies on a max-plus variational formulation, and exploits the properties of projectors on max-plus semimodules. We prove here a convergence result, in arbitrary dimension, showing that for a subclass of problems, the error estimate is of order $\delta+\Delta x(\delta)^{-1}$, where $\delta$ and $\Delta x$ are the time and space steps respectively. We also show how the max-plus analogues of the mass and stiffness matrices can be computed by convex optimization, even when the global problem is non convex. We illustrate the method by numerical examples in dimension 2.
\end{abstract}
\maketitle
\section{Introduction}
We consider the optimal control problem:
\begin{subequations}
\label{problemP}
\begin{align}
\label{p1}
\mrm{maximize }
\int_0^T \ell(x(s),u(s))\dd s+\phi(x(T))
\end{align}
over the set of trajectories $(x(\cdot),u(\cdot))$ satisfying
\begin{align}\label{p2}
\begin{split}
\dot{x}(s)&=f(x(s),u(s)),\quad
x(0)=x,\\
& x(s)\in X,\quad u(s)\in U,
\end{split}
\end{align}
\end{subequations}
for all $0\leq s\leq T$.
Here, the \new{state space} $X$ is a subset of $\R^n$, the set of \new{control values} $U$ is a subset of $\R^m$, the \new{horizon} $T>0$ and the \new{initial condition} $x\in X$ are given, we assume that the map $u(\cdot)$ is measurable, and that the map $x(\cdot)$ is absolutely continuous.
We also assume that the \new{instantaneous reward} or \new{Lagrangian} $\ell:X\times U \to \R$, and the \new{dynamics} $f:X\times U \to \R^n$, are sufficiently regular maps, and that the \new{terminal reward} $\phi$ is a map $X\to \R\cup\{-\infty\}$.

We are interested in the numerical computation of the \new{value function} $v$ which associates to any $(x,t)\in X\times [0,T]$ the supremum $v(x,t)$ of $\int_0^t \ell(x(s),u(s))\dd s+\phi(x(t))$, under the constraint~(\ref{p2}), for $0\leq s \leq t$.
It is known that, under certain regularity assumptions, $v$ is solution of the Hamilton-Jacobi equation
\begin{subequations}\label{HJ}
\begin{gather}
-\frac{\partial v}{\partial t}+H(x,\frac{\partial v}{\partial x})=0, \quad (x,t)
\in X \times (0,T] \enspace,\label{HJ1}\end{gather}
with initial condition:
\begin{gather}
v(x,0)=\phi(x), \quad  x \in X \enspace,
\end{gather}
\end{subequations}
where $H(x,p)=\sup_{u\in U}\ell(x,u)+p\cdot f(x,u)$ is the \new{Hamiltonian} of the problem (see for instance~\cite{lions,barles}).
The \new{evolution semigroup} $S^t$ of~(\ref{HJ}), or Lax-Oleinik semigroup, associates to any map $\phi$ the function $v^t:=v(\cdot,t)$, where $v$ is the value function of the optimal control problem~(\ref{problemP}). 

Maslov observed in~\cite{maslov73} that the evolution semigroup $S^t$ is max-plus linear (see also~\cite{maslov92,kolokoltsov}).
Recall that the \new{max-plus semiring}, $\rmax$, is the set $\R\cup\{-\infty\}$, equipped with the addition $a\oplus b=\max(a,b)$ and the multiplication $a\otimes b=a+b$. By \new{max-plus linearity}, we mean that for all maps $f,g$ from $X$ to $\rmax$, and for all $\lambda\in\rmax$, we have 
\begin{align*}
S^t(f\oplus g)&=S^tf\oplus S^tg \enspace ,\\
S^t(\lambda f)&=\lambda (S^tf) \enspace ,
\end{align*}
where $f\oplus g$ denotes the map $x\mapsto f(x)\oplus g(x)$, and $\lambda f$ denotes the map $x\mapsto \lambda \otimes f(x)$.
Linear operators over max-plus type semirings have been widely studied, see for instance~\cite{cuning,maslov92,baccelli,kolokoltsov,gondran-minoux}, see also~\cite{fathi}.

In~\cite{mceneaney}, Fleming and McEneaney introduced a first discretization method exploiting the max-plus linearity of the semigroup $S^t$.

In~\cite{MTNS}, we introduced a new max-plus based discretization method, inspired by the classical finite element method. The max-plus finite element method of~\cite{MTNS} approximates the evolution semigroup $S^t$ by means of a nonlinear discrete semigroup, which can be interpreted as the dynamic programming operator of a deterministic zero-sum two players game, with finite action and state spaces (unlike the method of Fleming and McEneaney which leads to a discrete optimal control problem). The state of the game corresponds to the set of finite elements. To each test function corresponds one possible action of the first player, and to each finite element corresponds one possible action of the second player. This discretization, which can be interpreted geometrically in terms of projections on semimodules, is similar to the classical Petrov-Galerkin finite element method.
 
The computation of the instantaneous payments of the game requires the evaluation of the max-plus scalar product $\<z,S^\delta w>$ for each finite element $w$ and each test function $z$, where $\delta$ is the time discretization step. In some special cases, $\<z,S^\delta w>$ can be computed analytically. In general, we need to approximate this scalar product, for each finite element $w$ and test function $z$. In~\cite{MTNS}, we used the simplest approximation $S^{\delta}w=w+\delta H(\cdot,\frac{\partial w}{\partial x})$, already considered in~\cite{mceneaney99}. This requires regularity assumptions on $w$ (or alternatively, on $z$, if one uses the dual semigroup~\cite{MTNS}). In this paper, we rather use a direct method, which allows us to approximate $\<z,S^\delta w>$ by the value of an optimization problem in finite dimension. We show that, under reasonable assumptions on $\ell$, $f$, $z$ and $w$, this approximation leads to a concave optimization problem. We also give an error estimate of order $\delta+\frac{\Delta x}{\delta}$.

The paper is organised as follows. In Section~\ref{prelim}, we recall some basic tools and notions: residuation, semimodules and projection. In Section~\ref{method}, we recall the formulation of the max-plus finite element method. The contents of Sections~\ref{prelim} and~\ref{method} are essentially taken from~\cite{MTNS}: we need to recall them to state our results. In Section~\ref{concave}, we discuss the approximation of the scalar product $\<z,S^\delta w>$. In Section~\ref{erreur}, we give the main convergence theorem. Finally, in Section~\ref{Numerical results}, we illustrate the method by numerical examples in dimension $2$.
\section{Preliminaries on residuation and projections over semimodules}\label{prelim}
In this section we recall some classical residuation results (see for example \cite{blyth}, \cite{baccelli}), and their application to linear maps on idempotent semimodules (see~\cite{litvinov,ilade}). We also review some results of \cite{wodes,ilade} concerning projectors over semimodules.
\subsection{Residuation, semimodules, and linear maps}
If $(S,\leq)$ and $(T,\leq)$ are (partially) ordered sets, we say that a map $f:S\to T$ is \new{monotone} if $s \leq s' \implies f(s)\leq f(s')$. We say that $f$ is \new{residuated} if there exists a map $f\sh: T\to S$ such that
\[
f(s) \leq t \iff s\leq f\sh(t) \enspace .
\]
The map $f$ is residuated if, and only if, for all $t\in T$, $\set{s\in S}{f(s)\leq t}$ has a maximum element in $S$. Then, 
\begin{align*}
f\sh(t)&=\max\set{s\in S}{f(s)\leq t},\quad \forall t\in T \enspace.
\end{align*}
If a set $\sK$ is a monoid for a commutative idempotent law $\oplus$ (\new{idempotent} means that $a\oplus a =a$), the \new{natural order} on $\sK$ is defined by $a\leq b \iff a\oplus b=b$. We say that $\sK$ is \new{complete} as a naturally ordered set if any subset of $\sK$ has a least upper bound for the natural order. If $(\sK,\oplus,\otimes)$ is an idempotent semiring, i.e., a semiring whose addition is idempotent, we say that the semiring $\sK$ is \new{complete} if it is complete as a naturally ordered set, and if the left and right multiplications, $L^\sK_a$, $R^\sK_a: \sK\to \sK$, $L^\sK_a(x)=ax$, $R^\sK_a(x)=xa$, are residuated.

The max-plus semiring, $\rmax$, is an idempotent semiring. It is not complete, but it can be embedded in the complete idempotent semiring  $\rmaxb$ obtained by adjoining $+\infty$ to $\rmax$, with the convention that $-\infty$ is absorbing for the multiplication $a\otimes b=a+b$.
The map $x\mapsto -x$ from $\rbar$ to itself yields an isomorphism from $\rmaxb$ to the complete idempotent semiring $\rminb$, obtained by replacing $\max$ by $\min$ and by exchanging the roles of $+\infty$ and $-\infty$ in the definition of $\rmaxb$.

Semimodules over semirings are defined like modules over rings, mutatis mutandis, see~\cite{litvinov,ilade}.
When $\sK$ is a complete idempotent semiring, we say that a (right) $\sK$-semimodule $\calx$ is \new{complete} if it is complete as a naturally ordered set, and if, for all $u\in \calx$ and $\lambda\in \sK$, the right and left multiplications, $R^\calx_{\lambda}:\;\calx\to \calx$, $v\mapsto v\lambda$ and $L^\calx_{u}:\;\sK\to \calx$, $\mu\mapsto u\mu$, are residuated.
In a complete semimodule $\calx$, we define, for all $u,v\in \calx$,
\begin{align*}
  u\lres v &\bydef (L_u^\calx)\sh(v) = \max\set{\lambda\in \sK}{u\lambda \leq v} \enspace .
\end{align*}
We shall use \new{semimodules of functions}: when $X$ is a set and $(\sK,\oplus,\otimes)$ is a complete idempotent semiring, the set of functions $\sK^X$ is a complete $\sK$-semimodule for the componentwise addition $(u,v)\mapsto u\oplus v$ (defined by $(u\oplus v)(x)= u(x)\oplus v(x)$), and the componentwise multiplication $(\lambda,u)\mapsto u \lambda$ (defined by $(u\lambda)(x)= u(x)\otimes\lambda$).

If $\sK$ is an idempotent semiring, and if $\calx$ and $\caly$ are $\sK$-semimodules, we say that a map $A:\calx\to \caly$ is \new{additive} if for all $u,v\in \calx$, $A(u\oplus v)=A(u)\oplus A(v)$ and that $A$  is \new{homogeneous} if for all $u\in \calx$ and $\lambda\in \sK$, $A (u \lambda)= A(u)\lambda$.
We say that $A$ is \new{linear}, or is a \new{linear operator}, if it is additive and homogeneous.
Then, as in classical algebra, we use the notation $Au$ instead of $A(u)$. 
When $A$ is residuated and $v\in \caly$, we use the notation $A\backslash v$ or $A\sh v$ instead of $A\sh (v)$.

If $X$ and $Y$ are two sets, $(\sK,\oplus,\otimes)$ is a complete idempotent semiring, and $a\in \sK^{X\times Y}$, we construct the linear operator $A$ from $\sK^Y$ to $\sK^X$ which associates to any $u\in \sK^Y$ the function $Au\in \sK^X$ such that $Au(x)=\supp_{y\in Y} a(x,y)\otimes u(y)$, where $\vee$ denotes the supremum for the natural order.
We say that $A$ is the \new{kernel operator} with \new{kernel} or \new{matrix} $a$.
We shall often use the same notation $A$ for the operator and the kernel.
As is well known (see for instance~\cite{baccelli}), the kernel operator $A$ is residuated, and\[
(A\backslash v)(y)=\inff_{x\in X}A(x,y)\backslash v(x),
\]
where $\wedge$ denotes the infimum for the natural order.
In particular, when $\sK=\rmaxb$, we have
\begin{align}
\label{e-conv}
(A\backslash v)(y)=\inff_{x\in X}(-A(x,y)+ v(x))= [- A^* (-v)](y)
\end{align}
where $A^*$ denotes the \new{transposed operator} $\sK^X\to \sK^Y$, which is associated to the kernel $A^*(y,x)=A(x,y)$. 
(In~(\ref{e-conv}), we use the convention that $+\infty$ is absorbing for addition.)

\subsection{Projectors on semimodules}
Let $\calv$ denote a \new{complete subsemimodule} of a complete semimodule $\calx$ over a complete idempotent semiring $\sK$, i.e., a subset of $\calx$ that is stable by arbitrary sups and by the action of scalars.
We call \new{canonical projector} on $\calv$ the map 
\begin{equation}\label{projecteur} P_\calv: \calx\to \calx, \quad u\mapsto P_\calv(u) = \maxx\set{v\in \calv}{v\leq u}. 
\end{equation}
Let $W$ denote a \new{generating family} of a complete subsemimodule $\calv$, which means that any element $v\in \calv$ can be written as $v=\supp\set{w\lambda_w}{w\in W}$, for some $\lambda_w\in\sK$.
It is known that
\[
P_\calv(u) =  \supp_{w\in W} w (w\lres u) 
\]
(see for instance~\cite{ilade}).
If $B:\calu\to\calx$ is a residuated linear operator, then the image $\im B$ of $B$ is a complete subsemimodule of $\calx$, and 
\begin{equation}\label{PimB}
P_{\im B}=B\comp B\sh.
\end{equation}
The max-plus finite element methods relies on the notion of projection on an image, parallel to a kernel, which was introduced by Cohen, the second author, and Quadrat, in~\cite{wodes}.
The following theorem, of which Proposition~\ref{vhdelta} below is an immediate corollary, is a variation on the results of~\cite[Section~6]{wodes}.
\begin{thm}[Projection on an image parallel to a kernel]\label{piBC}
Let $B:\calu\to\calx$ and $C:\calx\to\caly$ be two residuated linear operators. Let $\projimker{B}{C}=B\comp(C\comp B)\sh\comp C$. We have $\projimker{B}{C}=\projimker{B}{}\comp \projimker{}{C}$, where $\projimker{B}{}=B\comp B\sh$ and $\projimker{}{C}=C\sh\comp C$. Moreover, $\projimker{B}{C}$ is a projector $\big((\projimker{B}{C})^2=\projimker{B}{C}\big)$, and for all $x\in\calx$:
\[
\projimker{B}{C}(x)=\maxx\set{y\in\im B}{Cy\leq Cx}.
\]
\end{thm}
The results of~\cite{wodes} characterize the existence and uniqueness, for all $x\in X$, of $y\in \im B$ such that $Cy=Cx$. In that case, $y=\projimker{B}{C}(x)$.

When $\sK=\rmaxb$, and $C:\rmaxb^X\to\rmaxb^Y$ is a kernel operator, $\projimker{}C=C\sh\comp C$ has an interpretation similar to~(\ref{PimB}):
\[ 
\projimker{}{C}(v)=C\sh\comp C(v)=-P_{\im C^*}(-v)=P_{-\im C^*} (v)\enspace,
\]
where $-\im C^*$ is thought of as a $\rminb$-subsemimodule of $\rminb^X$, so that, 
\[
P_{-\im C^*} (v)= \min\set{w\in -\im C^*}{w\geq v} \enspace .
\]
where $\leq$ denotes here the usual order on $\rbar^X$, since the natural order of $\rminb^X$ is the reverse of the usual order.
When $B:\rmaxb^U\to\rmaxb^X$ is also a kernel operator, we have
\[
\projimker BC=P_{\im B}\comp P_{-\im C^*} \enspace .
\]
This factorization is instrumental in the geometrical interpretation of the finite element algorithm, see~\cite[Example 10]{MTNS}. 
\section{The max-plus finite element method}\label{method}
In this section we describe the max-plus finite element method to solve Problem~(\ref{problemP}).
Let $S^t$ and $v^t$ be defined as in the introduction. Since $S^{t+t'}=S^t\circ S^{t'}$, for $t,t'>0$, we obtain the recursive equation:
\begin{equation}\label{exact}
\begin{array}{cc}
v^{t+\delta}=S^{\delta} v^t, & t=0,\delta, \cdots, T-\delta
\end{array}
\end{equation}
with $v^0=\phi$ and $\delta=\frac{T}{N}$, for some positive integer $N$. Let $\calw$ be a $\rmaxb$-semimodule of functions from $X$ to $\rmaxb$ such that $\phi\in\calw$ and for all $v\in\calw$, $t>0$, $S^tv\in\calw$. We suppose given a ``dual'' semimodule $\calz$ of ``test functions'' from $X$ to $\rmaxb$. The max-plus \new{scalar product} is defined by $\<u,v>=\supp_{x\in X} u(x)\otimes v(x)$, for all functions $u,v:X\to \rmaxb$.
We replace~(\ref{exact}) by:
\begin{equation}\label{produitscalaire}
\<z,v^{t+\delta}>=\<z,S^{\delta}v^{t}> ,
\quad  \forall z\in\calz \enspace,
\end{equation}
for $t=0,\delta,\ldots,T-\delta$, with $v^{\delta},\ldots,v^T\in\calw$.
This equation can be seen as the analogue of a \new{variational} or \new{weak formulation}. Kolokoltsov and Maslov used this formulation in \cite{kolokltsovmaslov88} to define a notion of generalized solution of Hamilton-Jacobi equations.
We consider now a semimodule $\calw_h\subset\calw$ generated by the family $\{w_i\}_{1\leq i\leq p}$. We call \new{finite elements} the functions $w_i$. We approximate $v^t$ by $v_h^t\in\calw_h$, that is, $v^t\simeq v_h^t=\bigoplus_{i=1}^{p}w_i\lambda^t_i$, where $\lambda_i^t\in\rmax$. We also consider a semimodule $\calz_h\subset\calz$ generated by the family $\{z_j\}_{1\leq j\leq q}$. The functions $z_1,\cdots,z_q$ will act as \new{test functions}. We replace~(\ref{produitscalaire}) by
\begin{equation}\label{pdtscalapproch}
\begin{array}{cc}
\<z_j,v_h^{t+\delta}>=\<z_j,S^{\delta}v_h^{t}>, & \forall 1\leq j\leq q\enspace,\\
\end{array}
\end{equation}
for $t=0,\delta, \cdots, T-\delta$, with $v_h^0=\phi_h\simeq\phi$ and $v_h^t\in\calw_h$, $t=0,\delta,\cdots,T$.

Since Equation~(\ref{pdtscalapproch}) need not have a solution, we look for the maximal subsolution, i.e.\ the maximal solution $v_h^{t+\delta}\in \calw_h$ of 
\begin{subequations}\label{inegal}
\begin{align}
\<z_j,v_h^{t+\delta}> \quad \leq \quad \<z_j,S^{\delta}v_h^{t}> \quad \forall 1\leq j\leq q\enspace. \label{infouegal}
\end{align}
We also take for the approximate value function $v_h^0$ at time $0$ the maximal solution $v_h^0\in \calw_h$ of
\begin{align}
v_h^{0}\leq v^0\enspace.\label{inegvh0}
\end{align}
\end{subequations}
Let us denote by $W_h$ the max-plus linear operator from $\rmax^p$ to $\calw$ with matrix $W_h=\col(w_{i})_{1\leq i\leq p}$, and by $Z_h^*$ the max-plus linear operator from $\calw$ to $\rmaxb^q$ whose transposed matrix is $Z_h=\col(z_j)_{1\leq j\leq q}$.
This means that $W_h\lambda=\bigoplus_{i=1}^{p}w_i \lambda_i$ for all $\lambda=(\lambda_i)_{i=1,\ldots, p}\in \rmax^p$, and $(Z_h^* v)_j=\< z_j,v>$ for all $v\in\calw$ and $j=1,\ldots,q$.
\begin{prop}[\cite{MTNS}]\label{vhdelta}
The maximal solution $v_h^{t+\delta}\in \calw_h$ of~(\ref{inegal}) is given by $v_h^{t+\delta}=S_{h}^{\delta}v_h^t$, where
\[
S_{h}^{\delta}=\projimker{W_h}{Z_h^*}\comp S^{\delta}\enspace.
\]
\end{prop}
The following proposition provides a recursive equation verified by the vector of coordinates $\lambda^t$ and is proved in~\cite{MTNS}.
\begin{prop}[\cite{MTNS}]
Let $v_h^{t}\in \calw_h$ be the  maximal solution of~(\ref{inegal}), for $t=0,\delta,\ldots, T$. 
Then,  for every $t=0,\delta,\ldots, T$, there exists $\lambda^t\in\rmax^p$ such that $v_h^t=W_h\lambda^t$. Moreover, the maximal $\lambda^t$ satisfying these conditions verifies the recursive equation
\begin{equation}\label{lambdat+dt}
\lambda^{t}=(Z_h^*W_h)\backslash (Z_h^*S^{\delta}W_h\lambda^{t-\delta}) \enspace,
\end{equation}
with the initial condition $\lambda^{0}=W_h\backslash \phi$.
\end{prop}
For $1\leq i \leq p$ and $1\leq j\leq q$, we define:
\begin{align}
(A_h)_{ji}&=\<z_j,w_i> \label{matrixA}\\
(B_{h})_{ji}&=\<z_j,S^{\delta}w_i> \label{matrixB}
\end{align}
$A_h$ and $B_h$ are respectively the matrices of the max-plus linear operators $Z_h^*W_h$ and $Z_h^*S^{\delta}W_h$. Equation~(\ref{lambdat+dt}) may be written explicitly, for $1\leq i\leq p$, as 
\[
\lambda^t_i = \min_{1\leq j\leq q}\Big( -(A_{h})_{ji}+ \max_{1\leq k\leq p}\big( (B_h)_{jk} +\lambda^{t-\delta}_k\big) \Big) \enspace .
\]
This recursion may be interpreted as the dynamic programming equation of a deterministic zero-sum two players game, with finite action and state spaces.

The ideal max-plus finite element method can be summarized as follows:
\begin{enumerate}
\item Choose $\delta=\frac{T}{N}$ and the finite elements $(w_i)_{1\leq i\leq p}$ and $(z_j)_{1\leq j\leq q}$, 
\item Compute the matrix $A_h$ by~(\ref{matrixA}) and the matrix $B_{h}$ by~(\ref{matrixB}),
\item Compute $\lambda^{0}=W_h\backslash\phi$ and  $v_h^{0}=W_h\lambda^{0}$.
\item For $t=\delta, 2\delta,\ldots,T$, compute $\lambda^{t}=A_h\backslash (B_{h}\lambda^{t-\delta})$ and  $v_h^{t}=W_h\lambda^{t}$.
\end{enumerate}
Then, $v_h^t$ approximates the value function at time $t$, $v^t$.

Fleming and McEneaney proposed a max-plus based method~\cite{mceneaney}, which also uses the linear formulation~(\ref{exact}). They approximated the evolution semigroup $S^t$ by a max-plus linear semigroup acting on a finitely generated semimodule of functions. A comparison of this method with the ideal max-plus finite element method appears in~\cite{MTNS}.
\section{Small time approximation of the Lax-Oleinik semigroup}\label{concave}
Computing $A_h$ from~(\ref{matrixA}) is an optimization problem, whose objective function is concave for natural choices of finite elements and test functions (see Section~\ref{erreur} below). 
Evaluating every scalar product $\<z,S^\delta w>$ leads to a new optimal control problem, which is simpler to approximate than Problem~(\ref{problemP}), because the horizon is small, and the functions $z$ and $w$ have a regularizing effect.
In~\cite{MTNS}, we proposed to use the following approximation of $S^\delta w$ provided by the Hamilton-Jacobi equation~(\ref{HJ1}):
\begin{equation}\label{stilde}
S^{\delta}w(x)\simeq w(x)+\delta H(x,\frac{\partial{w}}{\partial x}), \quad \mrm{for all } x\in X .
\end{equation}
 In this paper, we use the approximation of $S^\delta w$ by the function $[S^{\delta}w]\appr $ such that, for all $x\in X$
\begin{equation}\label{stilde2}
[S^{\delta}w]\appr (x)=\sup_{u\in U}\Big\{\delta \ell(x,u)+w\big(x+\delta f(x,u)\big)\Big\} .
\end{equation}
Let $[S^{\delta}W_h]\appr $ denotes the max-plus linear operator from $\rmax^p$ to $\calw$ with matrix $[S^{\delta}W_h]\appr =\col([S^{\delta}w_i]\appr )_{1\leq i\leq p}$. 
The above approximation of  $S^{\delta}w$ yields an approximation of the matrix $B_h$ by the matrix $B_h\appr:= Z_h^* [S^{\delta}W_h]\appr $, whose entries are given, for $1\leq i \leq p$ and $1\leq j\leq q$, by:
\begin{equation}\label{e-convenient}
(B_h\appr )_{ji}\!=\hspace{-1ex}\sup_{x\in X,u\in U}\hspace{-1ex}\{z_j(x)+w_i\big(x+\delta f(x,u)\big)+\delta\ell(x,u)\}.
\end{equation}
The following proposition shows that under assumptions on $\ell$, $f$, $z_j$ and $w_i$, computing the approximation~\eqref{e-convenient} is a concave maximization problem. In this case, one can compute the entries of the matrix $B_h\appr$ using standard convex optimization algorithms.
\begin{prop}
Let $X$ be a convex set of $\R^n$ and let $U$ be a convex set of $\R^m$. Assume that $z:\R^n\to\R\cup\{-\infty\}$ is concave, $\ell\in\mathcal{C}^2(X\times U,\R)$, $w\in\mathcal{C}^2(\R^n,\R)$ and $f$ is affine. Let $F=\frac{\partial f}{\partial x}$, $G=\frac{\partial f}{\partial u}$ and let $\|\cdot\|$ denotes the Euclidean norm of operators. Assume that there exist $\alpha,\beta,C>0$ such that $-CI_{n}\leq\frac{\partial^2w}{\partial x^2}\leq-\beta I_{n}$, $\frac{\partial^2\ell}{\partial u^2}\leq-\alpha I_{m}$, $\|\frac{\partial^2\ell}{\partial x^2}\|\leq C$ and $\|\frac{\partial^2\ell}{\partial x\partial u}\|\leq C$, where $I_n$ is the $n\times n$ identity matrix. Then there exists a constant $\delta_0=\delta_0(\alpha,\beta,C,\|F\|,\|G\|)$ such that, for all $\delta\leq\delta_0$, the function $(x,u)\mapsto z(x)+w\big(x+\delta f(x,u)\big)+\delta\ell(x,u)$ is concave.
\end{prop}
\begin{proof}
Since $z$ is concave, it suffices to prove that the function $(x,u)\mapsto b(x,u)=w\big(x+\delta f(x,u)\big)+\delta\ell(x,u)$ is concave. Since $w$, $\ell$, $f$ are $\mathcal{C}^2$ and $X$ and $U$ are convex sets, we must show that $h^*\nabla^2b(x,u)h\leq0$ for all $h\in\R^{n+m}$, $x\in X$, $u\in U$, where $\nabla^2b(x,u)$ denotes the Hessian of $b$ at $(x,u)$.
Set $x'=x+\delta f(x,u)$ and $h=\begin{pmatrix}h_1 \\ h_2\end{pmatrix}$ where $h_1\in\R^n$ and $h_2\in\R^m$. Then
\begin{multline*}
h^*\nabla^2b(x,u)h=h_2^*\Big(\delta\frac{\partial^2\ell}{\partial u^2}(x)+ \delta^2 G^*\frac{\partial^2w}{\partial x^2}(x')G\Big)h_2+\\
2\delta h_1^*\Big(\frac{\partial^2\ell}{\partial x\partial u}(x)+(I_n+\delta F)^*\frac{\partial^2w}{\partial x^2}(x')G\Big)h_2\\
\quad +h_1^*\Big(\delta\frac{\partial^2\ell}{\partial x^2}(x)+(I_n+\delta F)^*\frac{\partial^2w}{\partial x^2}(x')(I_n+\delta F)\Big)h_1\\
\leq\Big(\delta C(1+2\|F\|)-\beta\Big)\|h_1\|^2-\delta\alpha\|h_2\|^2+\\
2\delta C\big(1+\|G\|+\delta\|F\|\|G\|\big)\|h_1\|\|h_2\|.
\end{multline*}
We recognize a quadratic form in the variables $\|h_1\|$ and $\|h_2\|$. This quadratic form will keep a constant sign (which is negative) for non zero values of $h$ if its discriminant is negative, i.e.\ if
\begin{multline*}
\delta^3C^2\|F\|^2\|G\|^2+2\delta^2C^2\|F\|\|G\|\big(1+\|G\|\big)+\\
\delta C\Big(C\big(1+\|G\|\big)^2+\alpha\big(1+2\|F\|\big)\Big)\leq\alpha\beta.
\end{multline*}
This is the case in particular if
\begin{multline*}
3\mathrm{max}\bigg(\delta^3C^2\|F\|^2\|G\|^2,2\delta^2C^2\|F\|\|G\|\big(1+\|G\|\big),\\
\delta C\Big(C\big(1+\|G\|\big)^2+\alpha\big(1+2\|F\|\big)\Big)\bigg)\leq\alpha\beta,
\end{multline*}
and so, we can take
\begin{equation*}
\begin{split}
\delta_0=\mathrm{min}\bigg(\sqrt[3]{\frac{\alpha\beta}{3C^2\|F\|^2\|G\|^2}},\sqrt{\frac{\alpha\beta}{6C^2\|F\|\|G\|\big(1+\|G\|\big)}},\\
\frac{\alpha\beta}{3C\Big(C\big(1+\|G\|\big)^2+\alpha\big(1+2\|F\|\big)\Big)}\bigg).
\end{split}
\end{equation*}
\end{proof}
\section{Error analysis}\label{erreur}
We first recall a general lemma showing that the error of the finite element method is controlled by the projection errors, $\|\projimker{W_h}{}v^t-v^t\|_{\infty}$ and $\|\projimker{}{Z_h^*}v^t-v^t\|_{\infty}$, and by the approximation error, $\|[S^{\delta}w_i]\appr -S^{\delta}w_i\|_{\infty}$.
\begin{lem}[\cite{MTNS}]\label{error}
For $t=0,\delta,\cdots,T$, let $v^t$ be the value function at time $t$, and $v_h^t$ be its approximation given by the max-plus finite element method, implemented with the approximation $B_h\appr$ of $B_h$, given by~(\ref{e-convenient}).
We have
\begin{gather*}
\|v_h^T- v^T\|_{\infty}\leq (1+\frac{T}{\delta})\Big(\max_{1\leq i\leq p}\|[S^{\delta}w_i]\appr -S^{\delta}w_i\|_{\infty}\\
\quad +\sup_{t=0,\delta,\ldots,T}(\|\projimker{}{Z_h^*}v^t-v^t\|_{\infty}+\|\projimker{W_h}{}v^t-v^t\|_{\infty})\Big)
\end{gather*}
\end{lem}
The proof of this lemma uses the fact that projectors over max-plus semimodules are non-expansive in the sup-norm.

To state an error estimate, we fix a norm $\|\cdot\|$ on $\R^n$ and we make the following assumptions:
\begin{itemize}
\item[$(H1)$] The semigroup preserves the set of $\frac{1}{c}$-semiconvex functions, for some $c>0$.
\item[$(H2)$] $f:X\times U\to \R^n$ is bounded and Lipschitz continuous with respect to $x$: $\exists L_f>0,M_f>0$ such that
\begin{align*}
\!\!\!\!\| f(x,u)-f(y,u)\| &\leq L_f\| x-y \|  &\!\forall x,y\in X,\forall u \in U\\
\!\!\!\!\| f(x,u)\|&\leq M_f, &\!\forall x\in X, u \in U.
\end{align*}
\item[$(H3)$] $\ell:X\times U\to \R$ is bounded and Lipschitz continuous with respect to $x$: $\exists L_\ell>0,M_\ell>0$ such that
\begin{align*}
| \ell(x,u)-\ell(y,u) |&\leq L_\ell\| x-y \| & \forall x,y \in X, u \in U,\\
| \ell(x,u) |&\leq M_l, & \forall  x\in X, u \in U.
\end{align*}
\item[$(H4)$] $\phi:X \to \R$ is bounded and Lipschitz continuous.
\end{itemize}
Recall that a function $f$ is \new{$\frac 1 c$-semiconvex} if $f(x)+\frac 1{2c}\|x\|_2^2$, where $\|\cdot\|_2$ is the standard Euclidean norm of $\R^n$, is convex. Spaces of semiconvex functions were already used by Fleming and McEneaney~\cite{mceneaney}.

We shall use the following finite elements.
\begin{defin}[Lipschitz finite elements]
We call \new{Lipschitz finite element} centered at point $\hat x\in X$, with constant $A>0$, the function $w(x)=-A\|x-\hat x\|_1$ where $\|x\|_1=\sum_{i=1}^{n}|x_i|$ is the $l^1$-norm of $\R^n$.
\end{defin}
The family of Lipschitz finite elements of constant $A$ generates, in the max-plus sense, the semimodule of Lipschitz continuous functions of Lipschitz constant $A$ with respect to $\|\cdot\|_1$.
\begin{defin}[Quadratic finite elements]
We call \new{quad\-ratic finite element} centered at point $\hat x\in X$, with Hessian $\frac 1 c>0$, the function $w(x)=-\frac{1}{2c}\|x-\hat x\|_2^2$.
\end{defin}
When $X=\R^n$, the family of quadratic finite elements with Hessian $\frac 1 c$ generates, in the max-plus sense, the semimodule of lower-semicontinuous $\frac 1 c$-semiconvex functions.
\begin{lem}\label{appr-St}
Let $X$ be a compact of $\R^n$. We make assumptions \mrm{(H2)} and \mrm{(H3)}. Assume that $w$ and its derivative are both Lipschitz continuous. Then there exists $K_1>0$ such that $\|[S^\delta w]\appr-S^\delta w\|_{\infty}\leq K_1\delta^2$.
\end{lem}
\begin{proof}
Denote by $M_{Dw}$ and $M_{D^2w}$ the Lipschitz constants with respect to norm $\|\cdot\|$ of $w$ and its derivative respectively. We first show that there exists $K>0$ such that $S^\delta w(x)-[S^\delta w]\appr(x)\geq -K\delta^2$. For all $x\in X$ we have 
\begin{eqnarray*}
(S^\delta w)(x)&\geq\sup\Big\{\int_0^{\delta}\ell(x(s),u)ds+w(x(\delta))\quad |u\in U,x(0)=x,\\
& \dot x(s)=f\big(x(s),u\big),s\in[0,\delta]\Big\}.
\end{eqnarray*}
In other words, we bounded $S^tw$ from below by considering only constant controls.\\
Since
\begin{align*}
\begin{split}
\Big|\int_0^\delta [\ell(x(s),u)-\ell(x,u)]ds\Big|&\leq L_\ell\int_0^\delta\|x(s)-x\|ds\\
&\leq L_{\ell}\int_0^{\delta}sM_fds\\
&\leq\frac{1}{2}L_{\ell}M_f\delta^2,
\end{split}
\end{align*}
we obtain
\begin{multline*}
(S^\delta w)(x)\geq-L_{\ell}M_f\frac{\delta^2}{2}+\sup\Big\{\delta \ell(x,u)+w(x(\delta))|u\in U,x(0)=x,\\
\dot x(s)=f\big(x(s),u\big),s\in[0,\delta]\Big\}.
\end{multline*}
Moreover
\begin{align*}
\begin{split}
\big|w(x(\delta))-w(x+\delta f(x,u))\big|&\leq M_{Dw}\|\int_0^\delta f(x(s),u)ds-\delta f(x,u)\|\\
&\leq M_{Dw}L_fM_f\frac{\delta^2}{2}\enspace .
\end{split}
\end{align*}
We deduce that
\begin{equation*}
(S^\delta w)(x)-[S^\delta w]\appr(x)\geq-\big(L_{\ell}+M_{Dw}L_f\big)M_f\frac{\delta^2}{2}\enspace .
\end{equation*}
We now prove an opposite inequality. For $x\in X$ we have
\begin{multline*}
(S^\delta w)(x)=\sup\Big\{\int_0^\delta\ell(x(s),u(s))ds+w\big(x+\int_0^\delta  f(x(s),u(s))ds\big)\quad |\\
\quad u(s)\in U,\dot x(s)=f(x(s),u(s)),x(0)=x,s\in[0,\delta]\Big\}.
\end{multline*}
By the same arguments as before, we show that
\begin{multline*}
(S^\delta w)(x)\leq\big(L_\ell+M_{Dw}L_f\big)M_f\frac{\delta^2}{2}+\\
\sup_{u(s)\in U,\forall s}\Big\{\int_0^\delta\ell(x,u(s))ds+w\big(x+\int_0^\delta f(x,u(s))ds\big)\Big\}.
\end{multline*}
Using the fact that
\begin{multline*}
w\big(x+\int_0^\delta f(x,u(s))ds\big)\leq w(x)+\nabla w(x)\cdot\Big(\int_0^\delta f(x,u(s))ds\Big)+M_{D^2w}M_f^2\frac{\delta^2}{2}\enspace ,
\end{multline*}
\begin{equation*}
\begin{split}
w\big(x+\delta f(x,u)\big)&\geq w(x)+\nabla w(x)\cdot\Big(\delta f(x,u)\Big)-M_{D^2w}M_f^2\frac{\delta^2}{2}
\end{split}
\end{equation*}
and
\begin{equation*}
\begin{split}
\int_0^\delta\big[\ell(x,u(s))+\nabla w(x)\cdot f(x,u(s))\big]ds\leq\sup_{u\in U}\delta \ell(x,u)+\nabla w(x)\cdot\big(\delta f(x,u)\big),
\end{split}
\end{equation*}
we deduce that
\begin{equation*}
S^\delta w(x)-[S^\delta w]\appr(x)\leq(L_\ell+M_{Dw}L_f+2M_{D^2w}M_f)M_f\frac{\delta^2}{2}.
\end{equation*}
\end{proof}
Using Lemma~\ref{error}, Lemma~\ref{appr-St} and explicit estimates of the projection errors appearing in Lemma~\ref{error}, along the lines of~\cite{asma}, we derive the following convergence result.
\begin{thm}\label{th-main}
Let $\Omega$ be an open convex set of $\R^n$, $X=\overline\Omega$ and $\hat X=X+\mathrm{B}(0,cL)$, where $L>0$, $c>0$. Suppose that $X$ and $\hat X$ have regular grids $\mathcal{T}_{\Delta x}$ and $\hat{\mathcal{T}}_{\Delta x}$ respectively of size $\Delta x$.
We make assumptions \mrm{(H1)}-\mrm{(H4)}, and assume that the value function at time $t$, $v^t$, is $L$-Lipschitz continuous with respect to $\|\cdot\|_1$ and $\frac{1}{c}$-semiconvex for all $t>0$, with the same constant $c$ as in \mrm{(H1)}.
Let us choose quadratic finite elements $w_i$ of Hessian $\frac 1 c$, centered at the points of $\hat{\mathcal{T}}_{\Delta x}$.
Let us choose, as test functions $z_j$, the Lipschitz finite elements with constant $A\geq L$, centered at the points of $\mathcal{T}_{\Delta x}$. 
For $t=0,\delta,\ldots, T$, let $v_h^t$ be the approximation of $v^t$ given by the max-plus finite element method implemented with the approximation $B_h\appr$ of $B_h$ given by~\eqref{e-convenient}.
Then, there exists a constant $K>0$ such that
\[
\|v_h^T-v^T\|_{\infty}\leq K(\delta+\frac{\Delta x}{\delta})\enspace .
\]
\end{thm}
\begin{rem}
A different approximation of $B_h$ relying on~\eqref{stilde} was used in~\cite{MTNS,asma}. It is easier to implement. In particular, it avoids the numerical solution of optimization problems, as soon as explicit formul{\ae} are available for the Hamiltonian $H$ and the point of maximum of $x\mapsto z_j(x)+w_i(x)$, which is frequently the case. However, it only leads to an error of order $\sqrt\delta+\frac{\Delta x}{\delta}$.
\end{rem}
\section{Numerical results}\label{Numerical results}
We now present some results obtained by the method discussed above.
\begin{exmp}[Linear Quadratic Problem]\label{ex-lq}
We consider the case where $U=\R^2$, $X=\R^2$, $\phi\equiv 0$,
\[
\ell(x,u)=-\frac{x_1^2+x_2^2}{2}-\frac{u_1^2+u_2^2}{2}\enspace\mrm{ and } \quad f(x,u)=u \enspace .
\]
We choose quadratic finite elements $w_i$ and $z_j$ of Hessian $\frac 1 c$.
We represent in Figure~\ref{quadratique} the solution given by our algorithm in the case where $T=5$, $\delta=0.5$, $\Delta x=0.05$, $c=0.1$. 
\begin{figure}[thpb]
\begin{center}
\includegraphics[scale=0.38]{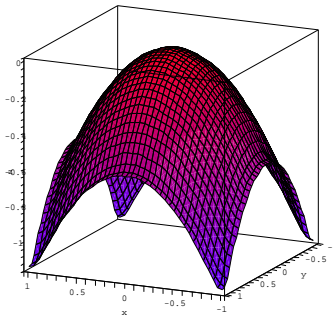}
\caption{Max-plus approximation of a linear quadratic control problem (Example~\ref{ex-lq})}
\label{quadratique}
\end{center}
\end{figure}

We observe a truncation effect on the boundaries of the set $X$. If we restrict $X$ to the set $[-0.5,0.5]^2$, we obtain a $L_{\infty}$-error of order 0.07.
\end{exmp}
\begin{exmp}[Distance problem]\label{ex-dist}
We consider the case where $T=1$, $\phi\equiv 0$, $X=[-1,1]^2$, $U=[-1,1]^2$,
\[
\ell(x,u)=\begin{cases}
-1 & \mathrm{if} \quad \|x\|_{\infty}< 1,\\
 0 & \mathrm{if} \quad \|x\|_{\infty}=1,
\end{cases}
\]
\[
f(x,u)=\begin{cases}
u & \mathrm{if} \quad \|x\|_{\infty}< 1,\\
0 & \mathrm{if} \quad \|x\|_{\infty}=1.
\end{cases}
\]
We choose quadratic finite elements $w_i$ of Hessian $\frac{1}{c}$ and Lipschitz finite elements $z_j$ with constant $A$. We represent in Figure~\ref{dist} the solution given by our algorithm in the case where $T=1$, $\delta=0.1$, $\Delta x=0.1$, $A=3$ and $c=1$. The $L_{\infty}$-error is of order $0.15$.
\begin{figure}[thpb]
\begin{center}
\includegraphics[scale=0.39]{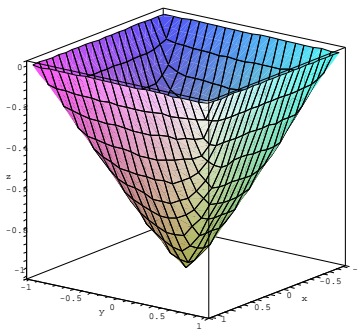}
\caption{Max-plus approximation of the distance problem (Example~\ref{ex-dist})}
\label{dist}
\end{center}
\end{figure}
\end{exmp}

\bibliographystyle{myalpha}

\end{document}